\documentclass[12pt,twoside]{amsart}
\usepackage{amssymb}
\usepackage{graphicx}
\usepackage[margin=1in]{geometry}

\DeclareMathOperator{\chara}{char}

\newtheorem{theorem}{Theorem}[section]
\newtheorem{lemma}[theorem]{Lemma}
\newtheorem{proposition}[theorem]{Proposition}
\newtheorem{corollary}[theorem]{Corollary}
\newtheorem{conjecture}[theorem]{Conjecture}

\theoremstyle{definition}
\newtheorem{remark}[theorem]{Remark}
\newtheorem{example}[theorem]{Example}

\usepackage{latexsym}
\usepackage{color,hyperref}
\definecolor{darkblue}{rgb}{0,0,0.6}
\hypersetup{colorlinks,breaklinks,linkcolor=darkblue,urlcolor=darkblue,anchorcolor=darkblue,citecolor=darkblue}

\title[Monomial Complete Intersections, WLP and Plane Partitions]{Monomial Complete Intersections, The Weak Lefschetz Property and Plane Partitions}

\author[\small Jizhou Li \tiny and \small Fabrizio Zanello]{\small Jizhou Li \tiny and \small Fabrizio Zanello}
\address{Department of Mathematical Sciences, Michigan Technological University, Houghton, MI 49931-1295}
\address{New address of the first author: Department of Computational and Applied Mathematics, Rice University, Houston, TX 77005-1892}
\email{jizhoul@mtu.edu, zanello@mtu.edu}

\begin{document}
\begin{abstract}
We characterize the  monomial complete intersections in three variables  satisfying the Weak Lefschetz Property (WLP), as a function of the characteristic of the base field. Our result  presents a surprising, and still combinatorially obscure, connection with the enumeration of plane partitions. It turns out that the rational primes $p$ dividing the number, $M \left( a,b,c \right)$, of plane partitions contained inside an arbitrary box of given  sides $a,b,c$ are precisely those for which a suitable  monomial complete intersection (explicitly constructed as a bijective function of $a,b,c$) fails to have the WLP in characteristic $p$. We wonder how powerful can be this connection between combinatorial commutative algebra and partition theory. We present a first result in this direction, by deducing, using  our algebraic techniques for the WLP,  some explicit information on the rational primes dividing $M \left( a,b,c \right)$.
\end{abstract}

\maketitle

\section{Introduction}

Let $ A=\bigoplus_{d\geq 0}A_{d} $ be a {\em standard graded  $K$-algebra}, where $K$ is an  infinite  field. $A$ can be identified with 
a quotient of a polynomial ring $K[x_1,\dots ,x_r]$ by a homogeneous ideal $I$, where the standard grading on $R$ (that is, all $x_i$'s have degree 1) is naturally induced on  $A$. The {\em Hilbert function} of $A$ is the arithmetic function $H$ defined by $H_A(d):=\dim_K A_d=\dim_KR_d-\dim_KI_d$, for all $d\geq 0$. We suppose here that $A=R/I$ be {\em artinian}. This has a number of equivalent formulations, including that the Krull-dimension of $A$ is zero, that the radical of $I$ is the irrelevant ideal $\underline{m}=(x_1,\dots ,x_r)$ of $R$, or that the Hilbert function of $A$ is eventually zero. This latter condition allows one to naturally identify $H_A$ with the {\em $h$-vector} of $A$, $h_A:=(h_0,h_1,...,h_e)$. Notice that $h_0=1$ and that we may assume, without loss of generality, that $h_e\neq 0$.

The {\em socle} of $A$ is the annihilator of $\underline{m}$ in $A$. Hence it is a homogeneous ideal,  and we define the {\em socle-vector} $s_A:=(s_0,...,s_e)$ to be its Hilbert function. It is easy to see that $s_e=h_e\neq 0$. The integer $e$ is defined as the {\em socle degree} of $A$ (or of $h_A$). If the socle is concentrated  in degree $e$, that is, $s_i=0$ for all $i\leq e-1$, we say that $A$ is a {\em level} algebra. If $A$ is level and $s_e=1$,  $A$ is called {\em Gorenstein}. (One often refers to the $h$-vector of a level or Gorenstein algebra as a level or Gorenstein $h$-vector.) The algebra $A$ is {\em monomial} if it is presented by  monomials (that is, if the ideal $I$ is generated by monomials).

Two  well-known and useful facts about Gorenstein algebras are that their $h$-vectors are  symmetric about the middle (i.e., $h_i=h_{e-i}$ for all indices $i$), and that if the algebras are also monomial and artinian, then they are  {\em complete intersections}. That is,  they are of the form $K[x_1, \dots ,x_r]/(x_1^{a_1}, \dots ,x_r^{a_r})$, for some positive integers $a_1, \dots ,a_r$.

One of the fundamental properties an artinian algebra can enjoy is the {\em Weak Lefschetz Property} (WLP). This is a very natural property, originally coming from algebraic geometry, which is also of great independent interest in algebra and combinatorics. $A$ is said to have
the WLP if there exists a linear form $L$ of $R$ such that, for all
indices $i$, the multiplication map $\times L$ between the $K$-vector spaces $A_{i}$ and $A_{i+1}$ has maximal rank. That is, $\times L$ is injective if $ \dim_K A_{i}\leq \dim_K A_{i+1} $ and surjective if $ \dim_K A_{i}\geq \dim_K A_{i+1} $ (and therefore bijective if $ \dim_K A_{i}=\dim_K A_{i+1} $). If such an $L$ exists, it is called a {\em Lefschetz element} of A. The Lefschetz elements of an algebra with the WLP form a non-empty open set in the Zariski topology of $\mathbb A^r(K)$, after we naturally identify a linear form with its coefficients.

A currently active and interesting line of research is to understand the behavior of the WLP for algebras over fields of positive characteristic (see \cite{3,5,ZZ}). Several of the initial results in this area have been unexpected or surprising --- especially in the light of what happens in characteristic zero --- and many problems today are still little understood. One of the main goals of this paper is to make a contribution in this direction.

We restrict to the case of monomial artinian Gorenstein quotients of a polynomial ring in $r=3$ variables, that is, algebras of the form $A=K[x,y,z]/(x^{\alpha},y^{\beta},z^{\gamma })$. Our main result  entirely characterizes  the positive integers $\alpha $, $\beta $, $\gamma $ and $p$, where $p$ is a  prime number, such that $A$ has the WLP in characteristic $p$. This  answers, as the particular case $\alpha =\beta =\gamma $, a question posed by Migliore, Mir\`o-Roig and Nagel (\cite{5}, Question 7.12). Also, for any such algebra $A$, the number of primes $p$ for which $A$ fails to have the WLP in characteristic $p$ is finite. In particular, this reproves a special case of a well-known result of Stanley \cite{St2}, saying that, in characteristic 0, all artinian monomial complete intersections have the WLP. (See also Watanabe \cite{watanabe} and Reid-Roberts-Roitman \cite{RRR}. Stanley's result actually showed much more, namely the {\em Strong} Lefschetz Property for such algebras.) At least one of the lemmas we prove along the way is also of some independent interest in terms of determinant evaluations.

As a byproduct, our main result yields a surprising connection with partition theory. It turns out that the  rational primes $p$ diving the number of plane partitions contained inside a given $a \times b \times c$ box  can be characterized as those prime numbers for which the monomial complete intersection $R=K[x,y,z]/(x^{a+b},y^{a+c},z^{b+c})$  fails to have the WLP in $\chara(K)=p$. It follows as the special case $a=1$ that the number, $\binom{b+c}{b}$, of integer partitions contained inside a $ b \times c$  rectangle is divisible by $p$ if and only if the algebra $K[x,y,z]/(x^{b+1},y^{c+1},z^{b+c})$  fails to have the WLP in $\chara(K)=p$. It would be very interesting to understand these facts also combinatorially.

We wonder how powerful  this new connection between combinatorial commutative algebra and partition theory can be for either field. We already move a first step in this direction, by  deducing, thanks to one of our algebraic techniques for the WLP, some  explicit information on the primes occurring in the integer factorization of  the number of plane partitions contained inside an arbitrary box.

\section{Preliminary results}

This section contains the preliminary results needed in the rest of the paper. The first of these is  known (as the Desnanot-Jacobi adjoint matrix theorem), and gives us a tool to compute the determinant of a matrix by induction.

\begin{lemma}[\cite{Jacobi}, Section 3; \cite{4}, Proposition 10]\label{proposition2.3}
Let $U$ be an $n\times n$ matrix. Denote by $U_{i_{1},\, i_{2},...,\, i_{k}}^{j_{1},\, j_{2},...,\, j_{k}}$ the
submatrix of $U$ in which rows $i_{1},\, i_{2},...,\, i_{k}$ and columns
$j_{1},\, j_{2},...,\, j_{k}$ are omitted. We have:

\[ \det(U)\cdot \det(U_{1,\, n}^{1,\, n})=\det(U_{1}^{1})\cdot \det(U_{n}^{n})-\det(U_{n}^{1})\cdot \det(U_{1}^{n}). \]
\end{lemma}

We now use  Lemma \ref{proposition2.3} to obtain the determinant of a particular matrix in closed form.

\begin{lemma}\label{lemma2.2}
 Let $N=\left(\binom{a+b}{a-i+j}\right)$, where $1\leq i\leq n+1,\,1\leq j\le n$.
Then, for any integer $k=1,\dots , n+1$, the determinant of the matrix $N_{k}=\left(\binom{a+b}{a-i+j}\right)$,
where $1\leq i\leq n+1,\,1\leq j\le n\, and\  i\not=k$, is
$$\det (N_{k})=\prod_{i=1}^{k-1}\tfrac{(n+1-i)(b+i)}{i(n+a-i)}\prod_{i=1}^{n}\tfrac{(a+b+i-1)!(i-1)!}{(a-2+i)!(b+i)!}.$$
 (As usual, we set any empty product to equal 1.)
\end{lemma}

\begin{proof} When $n=1$, it is clear that $N=\left[\begin{array}{c}
\binom{a+b}{b}\\
\binom{a+b}{b+1}\end{array}\right]$. If $k=1$, $\det(N_{k})=\binom{a+b}{b+1}$. If $k=2$,  $\det(N_{k})=\tfrac{b+1}{a}\cdot\tfrac{(a+b)!}{(a-1)!(b+1)!}=\binom{a+b}{b}$.\\

When $n=2$,
$$N=\left[\begin{array}{cc}
\binom{a+b}{b} & \binom{a+b}{b-1}\\
\binom{a+b}{b+1} & \binom{a+b}{b}\\
\binom{a+b}{b+2} & \binom{a+b}{b+1}\end{array}\right], N_{1}=\left[\begin{array}{cc}
\binom{a+b}{b+1} & \binom{a+b}{b}\\
\binom{a+b}{b+2} & \binom{a+b}{b+1}\end{array}\right], N_{2}=\left[\begin{array}{cc}
\binom{a+b}{b} & \binom{a+b}{b-1}\\
\binom{a+b}{b+2} & \binom{a+b}{b+1}\end{array}\right], N_{3}=\left[\begin{array}{cc}
\binom{a+b}{b} & \binom{a+b}{b-1}\\
\binom{a+b}{b+1} & \binom{a+b}{b}\end{array}\right].$$

One can easily check that
$$\det(N_{1})=\prod_{i=1}^{0}\tfrac{(n+1-i)(b+i)}{i(n+a-i)}\prod_{i=1}^{2}\tfrac{(a+b+i-1)!(i-1)!}{(a-2+i)!(b+i)!}, \det(N_{2})=\prod_{i=1}^{1}\tfrac{(n+1-i)(b+i)}{i(n+a-i)}\prod_{i=1}^{2}\tfrac{(a+b+i-1)!(i-1)!}{(a-2+i)!(b+i)!},$$ $$\det(N_{3})=\prod_{i=1}^{2}\tfrac{(n+1-i)(b+i)}{i(n+a-i)}\prod_{i=1}^{2}\tfrac{(a+b+i-1)!(i-1)!}{(a-2+i)!(b+i)!}.$$

We can now assume by induction that  the determinant be the desired one for square matrices $N_k$ of size up to $n-1$, for any given integer $n-1 \geq 2$. We want to show the result is true when $N_k$ has size $n$. We have:\\

$N_{k,1,n}^{1,n}=\left(\binom{a+b}{b-j+i}\right)$, which  is an $(n-2)\times(n-2)$ matrix with $i\not=k-1$;\\

$N_{k,1}^{1}=\left(\binom{a+b}{b-j+i}\right)$, which  is an $(n-1)\times(n-1)$ matrix with $i\not=k-1$; \\

$N_{k,n}^{n}=\left(\binom{a+b}{b-j+i}\right)$, which  is an $(n-1)\times(n-1)$ matrix with $i\not=k$;\\

$N_{k,1}^{n}=\left(\binom{a+b}{b+1-j+i}\right)$, which  is an $(n-1)\times(n-1)$ matrix with $i\not=k-1$;\\

$N_{k,n}^{1}=\left(\binom{a+b}{b-1-j+i}\right)$, which  is an $(n-1)\times(n-1)$ matrix with $i\not=k$.\\
\\\noindent
Thus, by induction, $\det(N_{k,1,n}^{1,n})=\prod_{i=1}^{k-2}\tfrac{(n-1-i)(b+i)}{i(n-2+a-i)}\prod_{i=1}^{n-2}\tfrac{(a+b+i-1)!(i-1)!}{(a-2+i)!(b+i)!}$;\\

$\det(N_{k,1}^{1})=\prod_{i=1}^{k-2}\tfrac{(n-i)(b+i)}{i(n-1+a-i)}\prod_{i=1}^{n-1}\tfrac{(a+b+i-1)!(i-1)!}{(a-2+i)!(b+i)!}$;\\

 $\det(N_{k,n}^{n})=\prod_{i=1}^{k-1}\tfrac{(n-i)(b+i)}{i(n-1+a-i)}\prod_{i=1}^{n-1}\tfrac{(a+b+i-1)!(i-1)!}{(a-2+i)!(b+i)!}$;\\

$\det(N_{k,1}^{n})=\prod_{i=1}^{k-2}\tfrac{(n-i)(b+1+i)}{i(n-2+a-i)}\prod_{i=1}^{n-1}\tfrac{(a+b+i-1)!(i-1)!}{(a-3+i)!(b+1+i)!}$;\\

$\det(N_{k,n}^{1})=\prod_{i=1}^{k-1}\tfrac{(n-i)(b+i-1)}{i(n+a-i)}\prod_{i=1}^{n-1}\tfrac{(a+b+i-1)!(i-1)!}{(a-1+i)!(b-1+i)!}$.\\
\\\noindent
A straightforward computation also gives:\\

{\smaller $$\tfrac{\det(N_{k,1}^{1})\cdot \det(N_{k,n}^{n})}{\det(N_{k,1,n}^{1,n})}=\tfrac{(b+n)(n+a-1)}{(1-k+n)(a+b+n-1)}\cdot\prod_{i=1}^{k-1}\tfrac{(n-i)(b+i)}{i(n+a-i)}\cdot\prod_{i=1}^{n}\tfrac{(a+b+i-1)!(i-1)!}{(a+i-2)!(b+i)!};$$}

{\smaller $$\tfrac{\det(N_{k,1}^{n})\cdot \det(N_{k,n}^{1})}{\det(N_{k,1,n}^{1,n})}=\tfrac{(-1+a)b}{(1-k+n)(a+b+n-1)}\cdot\prod_{i=1}^{k-1}\tfrac{(n-i)(b+i)}{i(n+a-i)}\cdot\prod_{i=1}^{n}\tfrac{(a+b+i-1)!(i-1)!}{(a-2+i)!(b+i)!}.$$}

Therefore,
$$\tfrac{\det(N_{k,1}^{1})\cdot \det(N_{k,n}^{n})}{\det(N_{k,1,n}^{1,n})}-\tfrac{\det(N_{k,1}^{n})\cdot \det(N_{k,n}^{1})}{\det(N_{k,1,n}^{1,n})}=\tfrac{n}{1-k+n}\prod_{i=1}^{k-1}\tfrac{(n-i)(b+i)}{i(n+a-i)}\prod_{i=1}^{n}\tfrac{(a+b+i-1)!(i-1)!}{(a+i-2)!(b+i)!}$$$$=\prod_{i=1}^{k-1}\tfrac{(n+1-i)(b+i)}{i(n+a-i)}\prod_{i=1}^{n}\tfrac{(a+b+i-1)!(i-1)!}{(a-2+i)!(b+i)!}=\det(N_{k}).$$

Hence, by Lemma \ref{proposition2.3}, $\det(N_{k})=\prod_{i=1}^{k-1}\tfrac{(n+1-i)(b+i)}{i(n+a-i)}\prod_{i=1}^{n}\tfrac{(a+b+i-1)!(i-1)!}{(a-2+i)!(b+i)!}.$
\end{proof}

\begin{remark}
Our previous lemma, which is also of  independent interest in terms of determinant evaluations, extends C. Krattenthaler's result on $ \underset{1 \leq i,j \leq n }{\det}\left(\binom{a+b}{a-i+j}\right) $ (see \cite{4}). The  determinant of the same matrix as Krattenthaler's had also been evaluated by P. Roberts \cite{robe} with a different method (thanks to Junzo Watanabe for kindly pointing out this reference). 
\end{remark}

The next lemma contains two very useful algebraic facts. In particular, it allows us to look at a unique map in order to determine whether our algebras $A$ have the WLP. (We state it here in a considerably lesser degree of generality than its  original formulation in \cite{5}.)

\begin{lemma}[\cite{5}, Propositions 2.1 and 2.2]\label{lemma2.4} Let $ A=R/I $ be an artinian monomial complete intersection  of socle degree $e$. Then $ L=x+y+z $ is a Lefschetz element of $A$. Also, $ A $ has the WLP if and only if the linear map $ \times L: A_{ \left\lfloor \tfrac{e-1}{2} \right\rfloor } \rightarrow  A_{ \left\lfloor \tfrac{e+1}{2} \right\rfloor } $ is injective.
\end{lemma}

Let us set $ s=\left\lfloor \tfrac{e-1}{2} \right\rfloor  $, and hence $ s+1=\left\lfloor \tfrac{e+1}{2} \right\rfloor $, for the rest of the paper. We call the integer $ s+1 $ a \textit{peak} of the $ h $-vector $h$, given that $h_{1}\leq h_{2} \leq\cdots\leq h_{s} \leq h_{s+1} \geq h_{s+2} \geq \cdots h_{e-1}\geq h_{e}=1$ (this chain of inequalities is well-known, and is essentially the {\em unimodality} property for complete intersection $h$-vectors). We prove next that, in the cases we are concerned with for $e$ even, $h$ has a single peak, i.e., $h_{s} < h_{s+1} > h_{s+2}$. When $e$ is odd,  by the symmetry of Gorenstein $h$-vectors, $h$ needs to have at least a twin peak, that is, $h_{s} = h_{s+1}$.

The following lemma tells us that, when it comes to determining when the WLP fails, we only need to be concerned with small values of $\gamma $ compared to $\alpha +\beta $. Precisely, we have:

\begin{lemma} Let $A=R/I$, where $R=K[x,y,z]$, $I=(x^{\alpha},y^{\beta},z^{\gamma})$, and $\alpha\leq\beta\leq\gamma $. Suppose that $\gamma >\alpha+\beta-2$ if $e$ is odd, and $\gamma >\alpha+\beta-3$ if $e$ is even. Then $A$  has the WLP.
\end{lemma}

\begin{proof}
It is easy to check that those values of $\gamma $ correspond to the case $\gamma > s+1$. Hence $A$ coincides through degree $s+1$ with the algebra $B=K[x,y,z]/(x^{\alpha},y^{\beta}) $. Since $z$ is clearly a non-zero divisor in $B$ (which has Krull-dimension 1), we have that multiplication by a general linear form  is an injective map between any two consecutive degrees of $B$. Hence, it is also injective through degree $s+1$ in $A$, which proves that $A$  has the WLP (cf. Lemma~\ref{lemma2.4}).
\end{proof}

\begin{remark}
The previous lemma takes care entirely of the case $\alpha =1$, that is, when our monomial complete intersections are essentially in two variables (after an obvious isomorphism). The stronger fact that \emph{all} algebras in two variables have the WLP in characteristic zero (and even more, the Strong Lefschetz Property) was first shown in \cite{HMNW}. The WLP part of that result was then reproved by J. Migliore and the second author \cite{MZ} using tools, including Green's theorem on hyperplane restrictions, which are independent of the characteristic (that the result was characteristic-free, however, was not mentioned in \cite{MZ}).
\end{remark}

\begin{lemma}\label{lemma2.6} Let $A=R/I$ be as above, and suppose that $\alpha\leq\beta\leq\gamma\leq\alpha+\beta-3$. If $e$ is even, then $h_A$ has a single peak at $s+1=\tfrac{\alpha+\beta+\gamma-3}{2}$, and $h_{s+1}-h_{s}=1$.
\end{lemma}

\begin{proof} Since for any $d$, $h_{d}=\binom{2+d}{d}-\dim_KI_d$, we have $$h_{s+1}-h_{s}=\left(\binom{s+3}{2}-\binom{s+2}{2}\right)-(\dim_KI_{s+1}-\dim_KI_s).$$
Hence it suffices to show that $\dim_KI_{s+1}-\dim_KI_s=s+1$.

The key observation here is that, from the assumption $\alpha\leq\beta\leq\gamma\leq\alpha+\beta-3$, it easily follows that no monomial in $I_{s+1}$ can be divisible by $x^{\alpha}y^{\beta}$ or $x^{\alpha}z^{\gamma}$ or $y^{\beta}z^{\gamma}$. Thus,

$$\dim_KI_{s+1}-\dim_KI_s=\left(\binom{s+1-\alpha +2}{2}+\binom{s+1-\beta +2}{2}+\binom{s+1-\gamma +2}{2}\right)-$$$$\left(\binom{s-\alpha +2}{2}+\binom{s-\beta +2}{2}+\binom{s-\gamma +2}{2}\right)=$$$$(s-\alpha +2)+(s-\beta +2)+(s-\gamma +2)=3s+3-(\alpha +\beta +\gamma -3)=3(s+1)-2(s+1)=s+1.$$
\end{proof}

\section{Monomial complete intersections in three variables}
\label{section3}

According to Lemma~\ref{lemma2.4}, determining when the WLP holds for $A=K[x,y,z]/(x^{\alpha},y^{\beta},z^{\gamma})$ is tantamount to determining  when the map $ \times L: A_{s} \rightarrow A_{s+1} $ is injective, for $ L=x+y+z $. We use linear algebra to study the problem.

First, we illustrate with a few examples the matrix $M$ associated with the map $\times L: A_{s} \rightarrow A_{s+1} $. Notice that $ M $ is a square matrix when  $h_A$ has a twin peak (that is, when the socle degree $e$ is odd), and by Lemma~\ref{lemma2.6}, it is of size $(h_{s}+1)\times h_{s}$ when $h_A$ has a single peak (i.e., for $e$ even). In the following examples the marked boxes inside $M$ are $1$'s and the rest of the matrix is $0$. They have been generated by means of a Mathematica \cite{6} computer program\footnote{The program generating the above matrices can be found on the first author's webpage, at\\
\url{http://www.mathlab.mtu.edu/~jizhoul/Commutative_Algebra_Project3.nb}}.\\

\begin{minipage}[h]{.23\linewidth}
\includegraphics[height=47mm,width=\linewidth]{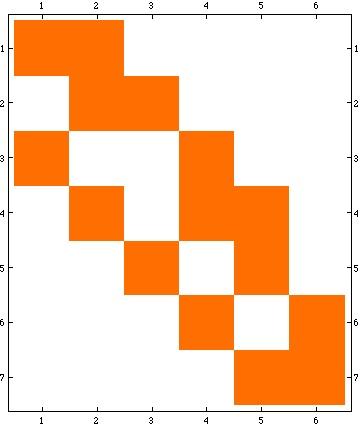}  \newline {\tiny Case $ \alpha=3, \beta=3, \gamma=3 $}
\end{minipage}
\begin{minipage}[h]{.23\linewidth}
\includegraphics[height=47mm,width=\linewidth]{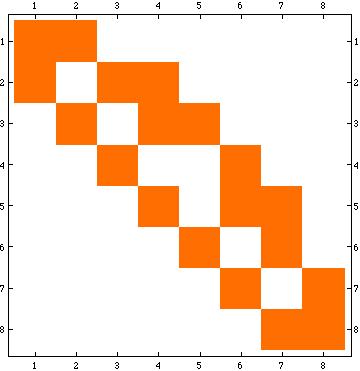} \newline {\tiny Case $ \alpha=3, \beta=3, \gamma=4 $}
\end{minipage}
\begin{minipage}[h]{.23\linewidth}
\includegraphics[height=47mm,width=\linewidth]{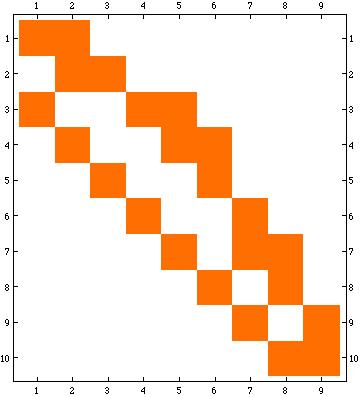} \newline {\tiny Case $ \alpha=3, \beta=4, \gamma=4 $}
\end{minipage}
\begin{minipage}[h]{.23\linewidth}
\includegraphics[height=47mm,width=\linewidth]{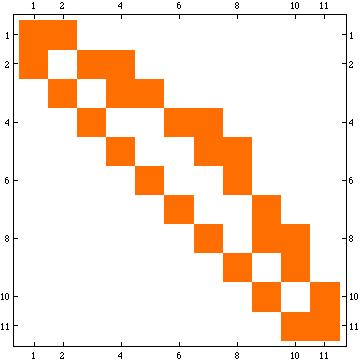} \newline {\tiny Case $ \alpha=3, \beta=4, \gamma=5 $}
\end{minipage}

\begin{minipage}[h]{.31\linewidth}
\includegraphics[height=60mm,width=\linewidth]{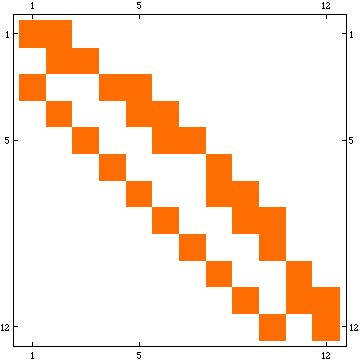} \newline {\tiny Case $ \alpha=4, \beta=4, \gamma=4 $}
\end{minipage}
\begin{minipage}[h]{.31\linewidth}
\includegraphics[height=60mm,width=\linewidth]{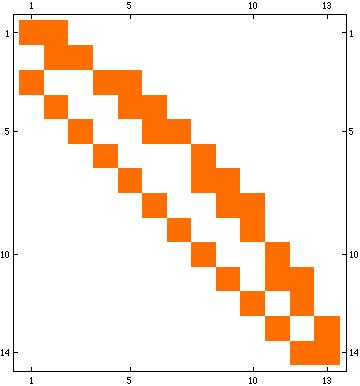} \newline {\tiny Case $ \alpha=4, \beta=4, \gamma=5 $}
\end{minipage}
\begin{minipage}[h]{.31\linewidth}
\includegraphics[height=60mm,width=\linewidth]{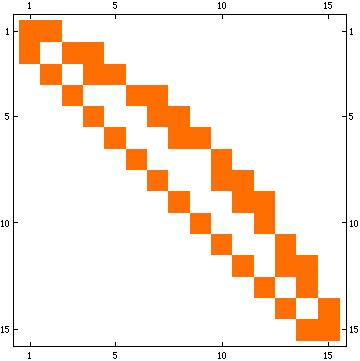} \newline {\tiny Case $ \alpha=4, \beta=4, \gamma=6 $}
\end{minipage}

\begin{minipage}[h]{.31\linewidth}
\includegraphics[height=60mm,width=\linewidth]{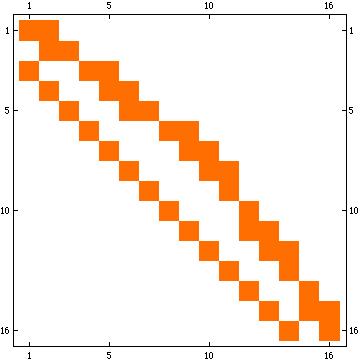} \newline {\tiny Case $ \alpha=4, \beta=5, \gamma=5 $}
\end{minipage}
\begin{minipage}[h]{.31\linewidth}
\includegraphics[height=60mm,width=\linewidth]{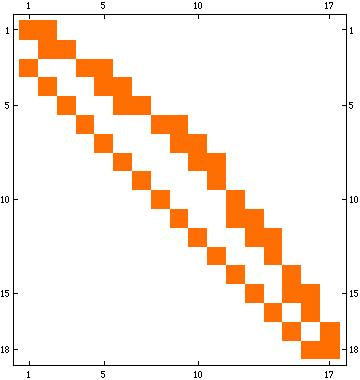} \newline {\tiny Case $ \alpha=4, \beta=5, \gamma=6 $}
\end{minipage}
\begin{minipage}[h]{.31\linewidth}
\includegraphics[height=60mm,width=\linewidth]{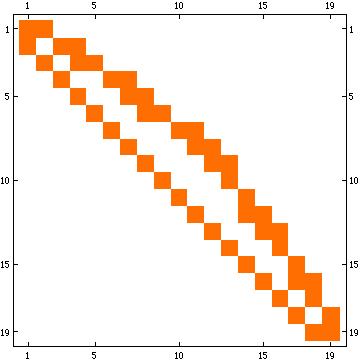} \newline {\tiny Case $ \alpha=4, \beta=5, \gamma=7 $}
\end{minipage}
{\ }\\
\\
\\

Now  let $$Z_{m\times(m+1)}=\left[\begin{array}{ccccccc}
1 & 1 & 0 & \cdots & 0 & 0 & 0\\
0 & 1 & 1 & \cdots & 0 & 0 & 0\\
\vdots & \vdots & \ddots & \ddots & \ddots & \vdots & \vdots\\
0 & 0 & 0 & \cdots & 1 & 1 & 0\\
0 & 0 & 0 & \cdots & 0 & 1 & 1\end{array}\right]_{m\times(m+1)}$$ 

and 
$$Z_{m\times m}=\left[\begin{array}{ccccccc}
1 & 1 & 0 & \cdots & 0 & 0 & 0\\
0 & 1 & 1 & \cdots & 0 & 0 & 0\\
\vdots & \vdots & \ddots & \ddots & \ddots & \vdots & \vdots\\
0 & 0 & 0 & \cdots & 1 & 1 & 0\\
0 & 0 & 0 & \cdots & 0 & 1 & 1\\
0 & 0 & 0 & \cdots & 0 & 0 & 1\end{array}\right]_{m\times m}.$$

We have the following result:

\begin{proposition} If $h_A$ has a twin peak, that is, $ \gamma = \alpha + \beta - 2m $ for some integer $m\geq 1$, then the matrix $M_{h_{s}\times h_{s}}$ is

{\tiny \[
\left[\begin{array}{ccccccccccc}
Z_{m\times(m+1)} & \textbf{0} & \cdots & \textbf{0} & \textbf{0} & \cdots & \textbf{0} & \textbf{0} & \cdots & \textbf{0} & \textbf{0}\\
I_{m+1} & Z_{(m+1)\times(m+2)} & \cdots & \textbf{0} & \textbf{0} & \cdots & \textbf{0} & \textbf{0} & \cdots & \textbf{0} & \textbf{0}\\
\vdots & \ddots & \ddots & \vdots & \vdots & \vdots & \vdots & \vdots & \cdots & \vdots & \vdots\\
\textbf{0} & \cdots & I_{\alpha-1} & Z_{(\alpha-1)\times\alpha} & \textbf{0} & \cdots & \textbf{0} & \textbf{0} & \cdots & \textbf{0} & \textbf{0}\\
\textbf{0} & \cdots & \cdots & I_{\alpha} & Z_{\alpha\times\alpha} & \cdots & \textbf{0} & \textbf{0} & \cdots & \textbf{0} & \textbf{0}\\
\vdots & \vdots & \vdots & \vdots & \vdots & \ddots & \vdots & \vdots & \cdots & \vdots & \vdots\\
\textbf{0} & \cdots & \cdots & \textbf{0} & \textbf{0} & \cdots & Z_{\alpha\times\alpha} & \textbf{0} & \cdots & \textbf{0} & \textbf{0}\\
\textbf{0} & \cdots & \cdots & \textbf{0} & \textbf{0} & \cdots & I_{\alpha} & Z_{(\alpha-1)\times\alpha}^{T} & \cdots & \textbf{0} & \textbf{0}\\
\vdots & \cdots & \cdots & \vdots & \vdots & \vdots & \vdots & \vdots & \ddots & \vdots & \vdots\\
\textbf{0} & \cdots & \cdots & \textbf{0} & \textbf{0} & \cdots & \textbf{0} & \textbf{0} & I_{m+2} & Z_{(m+1)\times(m+2)}^{T} & \textbf{0}\\
\textbf{0} & \cdots & \cdots & \textbf{0} & \textbf{0} & \cdots & \textbf{0} & \textbf{0} & \cdots & I_{m+1} & Z_{m\times(m+1)}^{T}\end{array}\right].\]
}
If $h_A$ has a single peak, that is, $ \gamma=\alpha + \beta -2m+1 $ for some integer $m> 1$, then $M_{(h_{s}+1) \times h_{s}}$ is 

{\tiny \[
\left[\begin{array}{ccccccccccc}
Z_{m\times(m+1)} & \textbf{0} & \cdots & \textbf{0} & \textbf{0} & \cdots & \textbf{0} & \textbf{0} & \cdots & \textbf{0} & \textbf{0}\\
I_{m+1} & Z_{(m+1)\times(m+2)} & \cdots & \textbf{0} & \textbf{0} & \cdots & \textbf{0} & \textbf{0} & \cdots & \textbf{0} & \textbf{0}\\
\vdots & \ddots & \ddots & \vdots & \vdots & \vdots & \vdots & \vdots & \cdots & \vdots & \vdots\\
\textbf{0} & \cdots & I_{\alpha-1} & Z_{(\alpha-1)\times\alpha} & \textbf{0} & \cdots & \textbf{0} & \textbf{0} & \cdots & \textbf{0} & \textbf{0}\\
\textbf{0} & \cdots & \cdots & I_{\alpha} & Z_{\alpha\times\alpha} & \cdots & \textbf{0} & \textbf{0} & \cdots & \textbf{0} & \textbf{0}\\
\vdots & \vdots & \vdots & \vdots & \vdots & \ddots & \vdots & \vdots & \cdots & \vdots & \vdots\\
\textbf{0} & \cdots & \cdots & \textbf{0} & \textbf{0} & \cdots & Z_{\alpha\times\alpha} & \textbf{0} & \cdots & \textbf{0} & \textbf{0}\\
\textbf{0} & \cdots & \cdots & \textbf{0} & \textbf{0} & \cdots & I_{\alpha} & Z_{(\alpha-1)\times\alpha}^{T} & \cdots & \textbf{0} & \textbf{0}\\
\vdots & \cdots & \cdots & \vdots & \vdots & \vdots & \vdots & \vdots & \cdots & \vdots & \vdots\\
\textbf{0} & \cdots & \cdots & \textbf{0} & \textbf{0} & \cdots & \textbf{0} & \textbf{0} & \ddots & \textbf{0} & \textbf{0}\\
\textbf{0} & \cdots & \cdots & \textbf{0} & \textbf{0} & \cdots & \textbf{0} & \textbf{0} & I_{m+1} & Z_{m\times(m+1)}^{T} & \textbf{0}\\
\textbf{0} & \cdots & \cdots & \textbf{0} & \textbf{0} & \cdots & \textbf{0} & \textbf{0} & \cdots & I_{m} & Z_{(m-1)\times m}^{T}\end{array}\right].\]
}
(There are $ \beta - \alpha $ blocks of $ Z_{\alpha \times \alpha} $ in both cases. $\textbf{0}$ represents a block matrix with all $0$'s.)
\end{proposition}

\begin{proof} We arrange the monomial basis of $ A_{s} $ in \textit{colexicographical order}. This ordering is defined, for monomials of the same degree, by setting $x_{1}^{a_{1}} x_{2}^{a_{2}} \cdots x_{n}^{a_{n}} < x_{1}^{a'_{1}} x_{2}^{a'_{2}} \cdots x_{n}^{a'_{n}}$ if  the following relation among base $10$ expansions is satisfied: $(a_{n} a_{n-1} \dots a_{1})_{10}< (a'_{n} a'_{n-1} \dots a'_{1})_{10}$.  For instance, the colexicographical order on the degree $3$ monomials of $K[x_1,x_2,x_3]$ gives $x_1^3<x_1^2 x_2<x_1 x_2^2<x_2^3<x_1^2 x_3<x_1x_2x_3<x_2^2 x_3<x_1 x_3^2<x_2 x_3^2<x_3^3$. 

We let $ \lbrace f_{s,1}, f_{s,2}, \dots , f_{s,h_{s}} \rbrace $ be the monomial basis of $ A_{s} $, arranged in colexicographical order. If $ f \in A_{s} $ then $ f=\alpha_{1} f_{s,1}+\alpha_{2} f_{s,2}+\dots +\alpha_{h_{s}} f_{s,h_{s}} $. By applying the linear operator $ \times L $ on $ f $, we get
$$ L\cdot f=(x+y+z) \cdot f = \beta_{1} f_{s+1,1}+\beta_{2} f_{s+1,2}+\dots +\beta_{h_{s+1}} f_{s+1,h_{s+1}} \in A_{s+1},$$
where $ \lbrace f_{s+1,1}, f_{s+1,2}, \dots , f_{s+1,h_{s+1}} \rbrace $ is the monomial basis of $ A_{s+1} $ in colexicographical order, and $ \beta_{i} = \alpha_{j}+\alpha_{k}+\alpha_{k+1} $. Notice that if $ x^{a}y^{b}z^{c} \in A_{s+1} $, then $ x^{a-1}y^{b}z^{c}, x^{a}y^{b-1}z^{c}, x^{a}y^{b}z^{c-1} \in A_{s} $, which implies that $ x^{a-1}y^{b}z^{c}, x^{a}y^{b-1}z^{c}$ are next to each other according to our ordering. Hence, the difference of the subscripts of the coefficients is $1$. It is easy to see that
$$M \cdot \left[\begin{array}{c}
\alpha_{1}\\
\alpha_{2}\\
\vdots\\
\alpha_{h_{s}}\end{array}\right]=\left[\begin{array}{c}
\beta_{1}\\
\beta_{2}\\
\vdots\\
\beta_{h_{s+1}}\end{array}\right].$$
If we now expand the product $ (x+y+z) \cdot f $ in terms of $ x, y, z$, after a standard but tedious computation we obtain that the matrix $ M $ has the desired form.
\end{proof}

Define, for any integers $a,b,c\geq 1$,

$$M \left(a,b,c \right):=\prod _{i=1}^a \prod _{j=1}^b \prod _{k=1}^c \tfrac{i+j+k-1}{i+j+k-2}=\prod _{i=1}^a \tfrac{(b+c+i-1)!(i-1)!}{(b+i-1)!(c+i-1)!},$$
and for $\alpha,\beta,\gamma $ as in our assumptions and such that $\alpha+\beta+\gamma $ is odd,
{\smaller $$ \mathcal{H}(k):=\prod_{i=1}^{k-1}\tfrac{\left(\tfrac{\alpha+\beta-\gamma+1}{2}-i\right)\left(\tfrac{-\alpha+\beta+\gamma-1}{2}+i\right)}{i(\alpha-i)}\cdot M \left( \tfrac{\alpha+\beta-\gamma-1}{2},\,\tfrac{\alpha-\beta+\gamma-1}{2},\tfrac{-\alpha+\beta+\gamma+1}{2} \right).$$}

We are now ready for the main result of this paper, where we characterize the artinian monomial complete intersections in $3$ variables having the WLP in characteristic $p$.

\begin{theorem}\label{main} Let $A=R/I$, where $R=K[x,y,z]$, $I=(x^{\alpha},y^{\beta},z^{\gamma})$,
and $\chara(K)=p$.\\
\begin{enumerate}
	\item If $e$ is odd (that is, $\alpha+\beta+\gamma $ is even), then $A$ fails to have 
the WLP if and only if $$p\mid M \left( \tfrac{\alpha+\beta-\gamma}{2},\,\tfrac{\alpha-\beta+\gamma}{2},\,\tfrac{-\alpha+\beta+\gamma}{2} \right).$$
\item If $e$ is even (that is, $\alpha+\beta+\gamma $ is odd), then $A$ fails to have the WLP if and only if, for all integers $1\leq k\leq \tfrac{\alpha+\beta-\gamma+1}{4}$,
$$p\mid\mathcal{H}(k).$$
\end{enumerate}
\end{theorem}

\begin{proof} Recall that we denote by $M$  the matrix associated with the map $\times L: A_{s} \rightarrow A_{s+1} $, where $L=x+y+z$. Our strategy  consists of computing the absolute value of the determinant of $M$ when $M$ is a square matrix, and all the maximal minors of $M$ when $M$ is not a square matrix. It follows from the above observations that $A$ fails to have the WLP in characteristic $p$ if and only if $p$ is a prime factor of the determinant of $M$ or of all the maximal minors of $M$. For simplicity we keep the notation $\det(T)$ when we actually mean the absolute value of $\det(T)$.

(1) Let $ \gamma=\alpha+\beta-2m $, for some $ m\geq 1 $.
Set {\tiny \[
U=\left[\begin{array}{cccccccccc}
I_{m+1} & Z_{(m+1)\times(m+2)} & \cdots & \textbf{0} & \textbf{0} & \cdots & \textbf{0} & \textbf{0} & \cdots & \textbf{0}\\
\vdots & \ddots & \ddots & \vdots & \vdots & \vdots & \vdots & \vdots & \cdots & \vdots\\
\textbf{0} & \cdots & I_{\alpha-1} & Z_{(\alpha-1)\times\alpha} & \textbf{0} & \cdots & \textbf{0} & \textbf{0} & \cdots & \textbf{0}\\
\textbf{0} & \cdots & \cdots & I_{\alpha} & Z_{\alpha\times\alpha} & \cdots & \textbf{0} & \textbf{0} & \cdots & \textbf{0}\\
\vdots & \vdots & \vdots & \vdots & \vdots & \ddots & \vdots & \vdots & \cdots & \vdots\\
\textbf{0} & \cdots & \cdots & \textbf{0} & \textbf{0} & \cdots & Z_{\alpha\times\alpha} & \textbf{0} & \cdots & \textbf{0}\\
\textbf{0} & \cdots & \cdots & \textbf{0} & \textbf{0} & \cdots & I_{\alpha} & Z_{(\alpha-1)\times\alpha}^{T} & \cdots & \textbf{0}\\
\vdots & \cdots & \cdots & \vdots & \vdots & \vdots & \vdots & \vdots & \ddots & \vdots\\
\textbf{0} & \cdots & \cdots & \textbf{0} & \textbf{0} & \cdots & \textbf{0} & \textbf{0} & I_{m+2} & Z_{(m+1)\times(m+2)}^{T}\\
\textbf{0} & \cdots & \cdots & \textbf{0} & \textbf{0} & \cdots & \textbf{0} & \textbf{0} & \cdots & I_{m+1}\end{array}\right],\]
}
{\tiny $
V=\left[\begin{array}{c}
\textbf{0}\\
\vdots\\
\textbf{0}\\
\textbf{0}\\
\vdots\\
\textbf{0}\\
\textbf{0}\\
\vdots\\
\textbf{0}\\
Z_{m\times(m+1)}^{T}\end{array}\right]$
},
{ $
W=\left[\begin{array}{cccccccccc}
Z_{m\times(m+1)} & \textbf{0} & \cdots & \textbf{0} & \textbf{0} & \cdots & \textbf{0} & \textbf{0} & \cdots & \textbf{0}\end{array}\right] $
},
{ $
X=\left[\textbf{0}\right] .$
}

Thus, {$M=\left[\begin{array}{cc}
W & X\\
U & V\end{array}\right]$}.
Evaluating $\det(M)$ up to sign is equivalent to evaluating the determinant of $M'=\left[\begin{array}{cc}
U & V\\
W & X\end{array}\right]$.
We have $ \det(M')=\det(U) \det(W U^{-1} V)=\det(W U^{-1} V) $. It is a standard task to check that $U^{-1}$ is the matrix

 {\tiny \[
\left[\begin{array}{cccc}
I_{m+1} & -Z_{(m+1)\times(m+2)} & \cdots & (-1)^{n}Z_{(m+1)\times(m+2)}Z_{(m+2)\times(m+3)}\cdots Z_{\alpha\times\alpha}\cdots Z_{\alpha\times\alpha}\cdots Z_{(m+2)\times(m+3)}^{T}Z_{(m+1)\times(m+2)}^{T}\\
\textbf{0} & I_{m+2} & \cdots & (-1)^{n-1}Z_{(m+2)\times(m+3)}\cdots Z_{\alpha\times\alpha}\cdots Z_{\alpha\times\alpha}\cdots Z_{(m+2)\times(m+3)}^{T}Z_{(m+1)\times(m+2)}^{T}\\
\vdots & \vdots & \ddots & \vdots\\
\textbf{0} & \textbf{0} & \textbf{0} & -Z_{(m+1)\times(m+2)}^{T}\\
\textbf{0} & \textbf{0} & \textbf{0} & I_{m+1}\end{array}\right].\]
}

Thus, {$ W U^{-1} V=(-1)^{n}Z_{m\times(m+1)}Z_{(m+1)\times(m+2)}\cdots Z_{\alpha\times\alpha}\cdots Z_{\alpha\times\alpha}\cdots Z_{(m+1)\times(m+2)}^{T}Z_{m\times(m+1)}^{T}.$} 
Using the properties of the Pascal triangle, one can easily verify that

 {\tiny $ Z_{m\times(m+1)}Z_{(m+1)\times(m+2)}\cdots Z_{\alpha\times\alpha}\cdots Z_{\alpha\times\alpha}\cdots Z_{(m+1)\times(m+2)}^{T}Z_{m\times(m+1)}^{T}=\left[\begin{array}{cccc}
\binom{\alpha+\beta-2m}{\beta-m} & \binom{\alpha+\beta-2m}{\beta-m-1} & \cdots & \binom{\alpha+\beta-2m}{\beta-2m+1}\\
\binom{\alpha+\beta-2m}{\beta-m+1} & \binom{\alpha+\beta-2m}{\beta-m} & \cdots & \binom{\alpha+\beta-2m}{\beta-2m+2}\\
\vdots & \vdots & \vdots & \vdots\\
\binom{\alpha+\beta-2m}{\beta-1} & \binom{\alpha+\beta-2m}{\beta-2} & \cdots & \binom{\alpha+\beta-2m}{\beta-m}\end{array}\right].$}

 By Krattenthaler's result (\cite{4}, 2.17), we have $$ \det \binom{\alpha+\beta-2m}{\beta-m+i-j}_{1 \leq i, j \leq m }=M \left( \tfrac{\alpha+\beta-\gamma}{2},\,\tfrac{\alpha-\beta+\gamma}{2},\,\tfrac{-\alpha+\beta+\gamma}{2} \right). $$
 
 Hence the absolute value of determinant of $ M $ is $M\left(\tfrac{\alpha+\beta-\gamma}{2},\,\tfrac{\alpha-\beta+\gamma}{2},\,\tfrac{-\alpha+\beta+\gamma}{2} \right),$ as desired.
 
(2) We now consider when $ M $ is a non-square matrix. In this case, we set $ \gamma=\alpha+\beta-2m+1 $, for some $ m>1 $.
 
 By Lemma~\ref{lemma2.6},  $ M $ is a $ (h_{s}+1) \times h_{s} $ matrix. Thus, its maximal minors are the matrices obtained by omitting one of the rows. Denote by $M_{k}$ the maximal minor obtained by omitting the $k$th row. We define $U$ to be\\
 
{\tiny $\left[\begin{array}{ccccccccccc}
I_{m+1} & Z_{(m+1)\times(m+2)} & \cdots & \textbf{0} & \textbf{0} & \cdots & \textbf{0} & \textbf{0} & \cdots & \textbf{0} & \textbf{0}\\
\vdots & \ddots & \ddots & \vdots & \vdots & \vdots & \vdots & \vdots & \cdots & \vdots & \vdots\\
\textbf{0} & \cdots & I_{\alpha-1} & Z_{(\alpha-1)\times\alpha} & \textbf{0} & \cdots & \textbf{0} & \textbf{0} & \cdots & \textbf{0} & \textbf{0}\\
\textbf{0} & \cdots & \cdots & I_{\alpha} & Z_{\alpha\times\alpha} & \cdots & \textbf{0} & \textbf{0} & \cdots & \textbf{0} & \textbf{0}\\
\vdots & \vdots & \vdots & \vdots & \vdots & \ddots & \vdots & \vdots & \cdots & \vdots & \vdots\\
\textbf{0} & \cdots & \cdots & \textbf{0} & \textbf{0} & \cdots & Z_{\alpha\times\alpha} & \textbf{0} & \cdots & \textbf{0} & \textbf{0}\\
\textbf{0} & \cdots & \cdots & \textbf{0} & \textbf{0} & \cdots & I_{\alpha} & Z_{(\alpha-1)\times\alpha}^{T} & \cdots & \textbf{0} & \textbf{0}\\
\vdots & \cdots & \cdots & \vdots & \vdots & \vdots & \vdots & \vdots & \ddots & \vdots & \vdots\\
\textbf{0} & \cdots & \cdots & \textbf{0} & \textbf{0} & \cdots & \textbf{0} & \textbf{0} & I_{m+2} & Z_{(m+1)\times(m+2)}^{T} & \textbf{0}\\
\textbf{0} & \cdots & \cdots & \textbf{0} & \textbf{0} & \cdots & \textbf{0} & \textbf{0} & \cdots & I_{m+1} & Z_{m\times(m+1)}^{T}\\
\textbf{0} & \cdots & \cdots & \textbf{0} & \textbf{0} & \cdots & \textbf{0} & \textbf{0} & \cdots & \textbf{0} & I_{m}\end{array}\right],$}
 
 {\tiny $V=\left[\begin{array}{c}
\textbf{0}\\
\vdots\\
\textbf{0}\\
Z_{(m-1)\times m}^{T}\end{array}\right]$}, { $W=\left[\begin{array}{ccccccccccc}
Z_{m\times(m+1)} & \textbf{0} & \cdots & \textbf{0} & \textbf{0} & \cdots & \textbf{0} & \textbf{0} & \cdots & \textbf{0} & \textbf{0}\end{array}\right]$}, { $X=\left[\textbf{0}\right]$}.

Hence {$M=\left[\begin{array}{cc}
W & X\\
U & V\end{array}\right]$}.
If $ 1 \leq k \leq m=\tfrac{\alpha+\beta-\gamma+1}{2}  $,  we need to omit the $k$th row of $W$ and $X$, in order to omit the $k$th row of $M$. Let $ W_{k} $ and $ X_{k} $ be the matrices $ W $ and $ X $ without their $k$th row. Therefore, $ M_{k} =\left[\begin{array}{cc}
W_{k} & X_{k}\\
U & V\end{array}\right]$. Obtaining $\det \left( M_{k} \right)  $ is equivalent to evaluating the determinant of $ \left[\begin{array}{cc}
U & V\\
W_{k} & X_{k}\end{array}\right]$.
Notice that $ \det \left[\begin{array}{cc}
U & V\\
W_{k} & X_{k}\end{array}\right]=\det (W_{k} U^{-1} V) $.\\

Hence, we can use the same approach as in part (1) in order to obtain the determinant of the matrix $ W_{k} U^{-1} V $. Entirely similar computations show that this matrix is {\footnotesize $$ W_{k} U^{-1} V=\left(\binom{\alpha+\beta -2m+1 }{ \beta - m+i-j}\right), \ \text{where} \ 1 \leq i \leq m,\ 1 \leq j \leq m-1,\  \text{and}\ i \neq k. $$} 

By Lemma~\ref{lemma2.2}, for $k=1,2, \dots \tfrac{\alpha+\beta-\gamma+1}{2}=m$, we have
$$ \det (W_{k} U^{-1} V)=\prod_{i=1}^{k-1}\tfrac{(\tfrac{\alpha+\beta-\gamma+1}{2}-i)(\tfrac{-\alpha+\beta+\gamma-1}{2}+i)}{i(\alpha-i)}\cdot M \left( \tfrac{\alpha+\beta-\gamma-1}{2},\,\tfrac{\alpha-\beta+\gamma-1}{2},\tfrac{-\alpha+\beta+\gamma+1}{2} \right)$$$$=\det( M_{k})=\mathcal{H}(k).$$

Notice that omitting the $k$th row, by symmetry, gives the same determinant (up to sign) that  we obtain by omitting the $\left(\tfrac{\alpha+\beta-\gamma+1}{2}+1-k\right)$th row.

However, we cannot apply the same method for $ k>m $, because in this case we would end up omitting a row from both matrices $ U $ and $ V $ and affecting the entire structure of $M_{k}$. We use a different approach to evaluate the determinant of $M_{k}$. Instead of deleting a row, we add an extra column $S$ with $1$ as the $k$th entry. The rest are all zeros. Now the new matrix $M'' $ is

\begin{flushleft}
{\tiny $\left[\begin{array}{cccccccccccc}
Z_{m\times(m+1)} & \textbf{0} & \cdots & \textbf{0} & \textbf{0} & \cdots & \textbf{0} & \textbf{0} & \cdots & \textbf{0} & \textbf{0} & \textbf{0}\\
I_{m+1} & Z_{(m+1)\times(m+2)} & \cdots & \textbf{0} & \textbf{0} & \cdots & \textbf{0} & \textbf{0} & \cdots & \textbf{0} & \textbf{0} & \textbf{0}\\
\vdots & \ddots & \ddots & \vdots & \vdots & \vdots & \vdots & \vdots & \cdots & \vdots & \vdots & S\\
\textbf{0} & \cdots & I_{\alpha-1} & Z_{(\alpha-1)\times\alpha} & \textbf{0} & \cdots & \textbf{0} & \textbf{0} & \cdots & \textbf{0} & \textbf{0} & \textbf{0}\\
\textbf{0} & \cdots & \cdots & I_{\alpha} & Z_{\alpha\times\alpha} & \cdots & \textbf{0} & \textbf{0} & \cdots & \textbf{0} & \textbf{0} & \textbf{0}\\
\vdots & \vdots & \vdots & \vdots & \vdots & \ddots & \vdots & \vdots & \cdots & \vdots & \vdots & \vdots\\
\textbf{0} & \cdots & \cdots & \textbf{0} & \textbf{0} & \cdots & Z_{\alpha\times\alpha} & \textbf{0} & \cdots & \textbf{0} & \textbf{0} & \textbf{0}\\
\textbf{0} & \cdots & \cdots & \textbf{0} & \textbf{0} & \cdots & I_{\alpha} & Z_{(\alpha-1)\times\alpha}^{T} & \cdots & \textbf{0} & \textbf{0} & \textbf{0}\\
\vdots & \cdots & \cdots & \vdots & \vdots & \vdots & \vdots & \vdots & \cdots & \vdots & \vdots & \vdots\\
\textbf{0} & \cdots & \cdots & \textbf{0} & \textbf{0} & \cdots & \textbf{0} & \textbf{0} & \ddots & \textbf{0} & \textbf{0} & \textbf{0}\\
\textbf{0} & \cdots & \cdots & \textbf{0} & \textbf{0} & \cdots & \textbf{0} & \textbf{0} & I_{m+1} & Z_{m\times(m+1)}^{T} & \textbf{0} & \textbf{0}\\
\textbf{0} & \cdots & \cdots & \textbf{0} & \textbf{0} & \cdots & \textbf{0} & \textbf{0} & \cdots & I_{m} & Z_{(m-1)\times m}^{T} & \textbf{0}\end{array}\right],$ }
\end{flushleft}

with
$$S=\left[\begin{array}{c}
0\\
\vdots\\
0\\
1\\
0\\
\vdots\\
0\end{array}\right].$$

Note that $ M'' $ is a square matrix. Evaluating  $\det( M_k) $ is equivalent to evaluating  $\det(M'') $. We again consider $ M'' $ as four blocks, where now
$$V=\left[\begin{array}{cc}
\textbf{0} & \textbf{0}\\
\vdots & S\\
\textbf{0} & \textbf{0}\\
\textbf{0} & \textbf{0}\\
\vdots & \vdots\\
\textbf{0} & \textbf{0}\\
\textbf{0} & \textbf{0}\\
\vdots & \vdots\\
\textbf{0} & \textbf{0}\\
\textbf{0} & \textbf{0}\\
Z_{(m-1)\times m}^{T} & \textbf{0}\end{array}\right],$$
and the other three blocks are the same as before. Therefore we only need to evaluate the determinant of $ W U^{-1} V $.

Employing the same method, after a series of standard computations we get

 $$\det( M_{k}) = \det (W U^{-1} V)=\sum_{k=1}^{\tfrac{\alpha+\beta-\gamma+1}{2}}(-1)^{k}\cdot\binom{s}{j-k}\cdot\mathcal{H}(k)$$
for $$\begin{cases} 1\leq s\leq\tfrac{\alpha-\beta+\gamma-1}{2}\ \ \ \ \ \ \ \ \ \ \ \ \ \text{and} & 1\leq j\leq\tfrac{\alpha+\beta-\gamma+1}{2}+s\\ \tfrac{\alpha-\beta+\gamma+1}{2}\leq s\leq\tfrac{-\alpha+\beta+\gamma-1}{2}\ \ \text{and}  & 1\leq j\leq\alpha\end{cases},$$
 \\
 and
 $$p\mid\sum_{k=1}^{\tfrac{\alpha+\beta-\gamma+1}{2}}(-1)^{k}\cdot\binom{s}{\tfrac{-\alpha+\beta+\gamma-1}{2}+k-j}\cdot\mathcal{H}(k)$$
 for all $\tfrac{-\alpha+\beta+\gamma+1}{2}\leq s\leq\gamma-1$ such that $1\leq j\leq\tfrac{\alpha+\beta+\gamma-1}{2}-s$, and for all $1\leq k\leq\tfrac{\alpha+\beta-\gamma+1}{2} $.
 
Notice that, if $p \mid \mathcal{H}(k)$, then obviously, for all $1\leq k\leq\tfrac{\alpha+\beta-\gamma+1}{2} $,
 $$p \mid \sum_{k=1}^{\tfrac{\alpha+\beta-\gamma+1}{2}}(-1)^{k}\cdot\binom{s}{j-k}\cdot\mathcal{H}(k)$$
and
 $$p \mid \sum_{k=1}^{\tfrac{\alpha+\beta-\gamma+1}{2}}(-1)^{k}\cdot\binom{s}{\tfrac{-\alpha+\beta+\gamma-1}{2}+k-j}\cdot\mathcal{H}(k).$$
 
Therefore $p$ only needs to satisfy the condition $p \mid \mathcal{H}(k)$ for all $1\leq k\leq\tfrac{\alpha+\beta-\gamma+1}{4} $. This completes  the proof of the theorem.
\end{proof}

\begin{remark}
In a different but similar effort to ours, M. Hara and J. Watanabe \cite{hawa} recently computed the determinants of the incidence matrices  between  graded components of the Boolean lattice on an $r$-set. This was equivalent, because of considerations analogous to those we have made above, to determining in which characteristics the complete intersections $K[x_1,\dots ,x_r]/(x_1^2,\dots ,x_r^2)$ have the Strong Lefschetz Property.
\end{remark}

\begin{example}
Let $ A=K[x,y,z]/(x^{b+1},y^{b+1},z^{2b}) $, where $ b$ is any positive integer. Hence the socle degree of $A$ is $4b-1$, and by Theorem \ref{main}, $A$ fails to have the WLP in characteristic $p$ if and only if $p$ divides $M \left( 1,b,b \right)=\binom{2b}{b}$.

Notice that $ 2\mid \binom{2b}{b} $ for all integers $b\geq 1$. For instance, this can be easily seen by observing that the involution on the class of $b$-subsets of a given $2b$-set, defined by taking the complementary of each set, has no fixed points. Thus, for any $b\geq 1$, $A$ fails to have the WLP in $ \chara(K)=2 $.
\end{example}

\begin{remark}
Notice that,  over the integers, the matrix $M$ of Theorem \ref{main} has a non-zero determinant when it is a square matrix. Similarly, when it is non-square, it is easy to see that not all its maximal minors $ M_{k} $ can have a zero determinant (for instance, $M_1$). This  reproves that all artinian monomial complete intersections in three variables have the WLP over a field of  characteristic zero, which is a (very) special case of a well-known result of Stanley \cite{St2} (see also Watanabe \cite{watanabe}, and for the first proof using only commutative algebra methods, Reid-Roberts-Roitman \cite{RRR}).
\end{remark}

Theorem \ref{main} also  answers, as the particular case $ \alpha=\beta=\gamma $, a question raised by Migliore,  Mir\`o-Roig and Nagel (\cite{5}, Question 7.12). Define, for any given $d\geq 3$ odd, a function $\mathcal{F}$ as
$$\mathcal{F}(k):=\prod_{i=1}^{k-1}\tfrac{\left(\tfrac{d+1}{2}-i\right)\cdot \left(\tfrac{d-1}{2}+i\right)}{i(d-i)}\cdot M \left( \tfrac{d-1}{2},\,\tfrac{d-1}{2},\,\tfrac{d+1}{2} \right).$$

\begin{corollary}\label{35} Let $A=R/I$, where $R=K[x,y,z]$, $I=(x^{d},y^{d},z^{d})$, and $\chara(K)=p$.
Then $A$ fails to have the WLP if and only if 
$p\mid M \left( \tfrac{d}{2},\ \tfrac{d}{2},\ \tfrac{d}{2} \right)$ if $d$ is even, and 
$p\mid\mathcal{F}(k)$ for all integers $1\leq k\leq \tfrac{d+1}{4}$ if  $d$ is odd.
\end{corollary}
\begin{proof}
This follows immediately from Theorem~\ref{main}, by setting $ \alpha=\beta=\gamma=d $.
\end{proof}

In general, even though Theorem \ref{main} has  established  necessary and sufficient conditions in order for an algebra to fail the WLP,  obviously determining explicitly such algebras is extremely difficult computationally. Indeed, this problem is equivalent to that of determining the primes in the integer factorization of the determinants of Theorem \ref{main}. However, we can prove with a simple algebraic argument the following explicit bounds in any characteristic $p$:

\begin{proposition}\label{tttt}
Let $A=K[x,y,z]/(x^{\alpha},y^{\beta},z^{\gamma })$, where $\alpha \leq \beta \leq \gamma $. Then, for any prime number $p$ such that
$$\gamma \leq p^n \leq  \left\lfloor \tfrac{\alpha +\beta +\gamma -3}{2} \right\rfloor $$
for some positive integer $n$, $A$ fails to have the WLP in $ \chara(K)=p $.
\end{proposition}

\begin{proof}
Recall that the peak of $h_A$ is  in degree $s+1= \left\lfloor \tfrac{\alpha +\beta +\gamma -3}{2} \right\rfloor$. It is a nice (combinatorial) exercise to check that $p \mid \binom{p^n}{k}$ for all integers $k=1,2, \dots , p^n -1$ (see \cite{bona}). Thus, from the assumption $\alpha \leq \beta \leq \gamma  \leq p^n $, we have that $$(x+y+z)\cdot (x+y+z)^{p^n-1}=(x+y+z)^{p^n}=0$$ in $A$. Therefore, since $L=x+y+z\neq 0$ in $A$, by induction on the degree we easily get that the map $ \times L: A_{p^n-1} \rightarrow A_{p^n} $ is not injective. Hence $A$ fails to have the WLP.
\end{proof}

\begin{example}
Consider again the special case $\alpha = \beta = \gamma =d$. The bounds of Proposition \ref{tttt} become
$$d\leq p^n\leq \left\lfloor \tfrac{3d -3}{2} \right\rfloor,$$
for some integer $n\geq 1$. It is easy to see that this already proves that at least one third of all algebras  $A=K[x,y,z]/(x^d,y^d,z^d) $ fail to have the WLP in a given $\chara(K)=p$.
\end{example}

We propose the following conjecture in characteristic $2$.

\begin{conjecture}\label{char2}
The algebra $ K[x,y,z]/(x^d,y^d,z^d) $  has the WLP in $ char(K)=2 $ if and only if $ d=\left\lfloor \tfrac{2^{n}+1}{3} \right\rfloor $ for some positive integer $n$.
\end{conjecture}

Notice, for instance, that $ d=\left\lfloor \tfrac{2^{n}+1}{3} \right\rfloor $ is odd for all $n\geq 1$, and therefore it  follows from our conjecture that $ K[x,y,z]/(x^d,y^d,z^d) $ has the WLP for all $d$ even. It would be very interesting to find a combinatorial proof of our conjecture, especially in the light of the connection we present in the next section between the WLP and the enumeration of plane partitions.

\section{The connection with plane partitions}\label{section4}

A {\em plane partition} of a positive integer $n$ is a finite two-dimensional array $A=(a_{i,j})$ of positive  integers, non-increasing from left to right and top to bottom, that add up to $n$. That is, $a_{i,j}\geq a_{i,j+1}\geq a_{i+1,j}\geq 1$ for all $i$ and $j$, and $\sum_{i,j}a_{i,j}=n$. (For details on this fascinating topic, see for instance \cite{And}.) 

We say that a plane partition $A=(a_{i,j})$ is {\em contained inside an $a \times b \times c$ box}, if $1\leq i\leq a$, $1\leq j\leq b$, and $a_{i,j}\leq c$ for all $i$ and $j$. P.A. MacMahon  determined the number of plane partitions contained inside an $a \times b \times c$ box (see \cite{Macdonald,Percy,pak}; for the first combinatorial proof of this result, see \cite{kra}). In fact, surprisingly, he proved that that number is $M \left( a,b,c \right)$, the very same $M \left( a,b,c \right)$  we   met in the previous section when determining the WLP  for our  monomial complete intersections.

\begin{remark}
A similar relationship, in that case involving certain classes of monomial almost complete intersections, has been discovered (independently but earlier) also by Cook and Nagel \cite{3}. (Their paper actually mentions other combinatorial objects, lozenges, which are known to be in bijection with plane partitions; see \cite{DT,DZZ}.) 
\end{remark}

Since $M \left( \tfrac{\alpha+\beta-\gamma}{2},\,\tfrac{\alpha-\beta+\gamma}{2},\,\tfrac{-\alpha+\beta+\gamma}{2} \right) $  enumerates, by MacMahon's result, the plane partitions contained inside an $ \tfrac{\alpha+\beta-\gamma}{2} \times \tfrac{\alpha-\beta+\gamma}{2} \times \tfrac{-\alpha+\beta+\gamma}{2} $ box, from Theorem \ref{main} we have:

\begin{theorem}\label{plus}
For any given  positive integers $ a,b, c $, the number of plane partitions contained inside an $a \times b \times c$  box is divisible by a rational prime $p$ if and only if the algebra $K[x,y,z]/(x^{a+b},y^{a+c},z^{b+c})$  fails to have the WLP when $\chara(K)=p$. 
\end{theorem}
\begin{proof}
This is immediate from Theorem~\ref{main}.  
\end{proof}

The case $a=1$ corresponds to that of  {\em integer partitions} (among the several possible choices, for an introduction,  a survey of the  main results and techniques, or  the philosophy behind this remarkably broad field, see \cite{And,AE,pak,St0}). Thus, from Theorem \ref{plus} we immediately have:

\begin{corollary}
For any given  positive integers $b$ and $ c $, the number of integer partitions, $M \left( 1,b,c \right)=\binom{b+c}{b}$, contained inside a{\ } $ b \times c$  rectangle is divisible by a rational prime $p$ if and only if the algebra $K[x,y,z]/(x^{b+1},y^{c+1},z^{b+c})$  fails to have the WLP when $\chara(K)=p$.
\end{corollary}

It seems reasonable to believe that such a nice connection between combinatorial commutative algebra and partition theory must have some deep combinatorial explanation. However, this is  unclear to us at the moment. Notice that our bijection is entirely different from the more natural one given by associating, to each monomial artinian ideal $I$ in three variables, the plane partition  whose solid Young diagram is the staircase diagram of (the order ideal of monomials outside of) $I$ (see Miller-Sturmfels \cite{MS} for details).

We wonder how powerful the connection between monomial complete intersections and plane partitions given by Theorem \ref{plus} could be for either field. One of our algebraic techniques used in studying the WLP allows us to move a first step in this direction, by providing a  highly non-trivial result on the function enumerating plane partitions. Namely, we are able to deduce some explicit information on the possible primes occurring in the integer factorization of  the number of plane partitions contained inside an arbitrary box. We have:

\begin{theorem}
{\ }\\
\begin{enumerate}\label{divide}
\item Fix three positive integers $a\leq b\leq c-1$. Then, for any prime number $p$ such that $$b+c\leq p^n\leq a+b+c-1$$ for some integer  $n\geq 1$, we have
$$p\mid M \left( a,b,c \right).$$
\item Fix two positive integers $a\leq b$. Then, for any prime number $p$ such that $$2b\leq p^n\leq a+2b-2$$ for some integer  $n\geq 1$, we have 
$$p\mid M \left( a,b,b \right).$$
\end{enumerate}
\end{theorem}

\begin{proof}
\begin{enumerate}
\item Set  $\alpha =a+b+1$, $\beta=a+c$ and $\gamma =b+c$. The result easily follows from Theorem \ref{main}, (2) (considering $\mathcal{H}(1)$), and Proposition \ref{tttt}.
\item Now set $\alpha =a+b$, $\beta=a+c$ and $\gamma =b+c$. The result  follows from Theorem \ref{main}, (1),  and Proposition \ref{tttt}.
\end{enumerate}
\end{proof}

\begin{example} Let $a=b=50$. Then Theorem \ref{divide} gives that the number, $M \left( 50,50,50 \right)$, of plane partitions contained inside a $50 \times 50 \times 50$ box is divisible by  any rational prime $p$ such that $100\leq p^n\leq 148$ for some $n\geq 1$. These values of $p$ are $$2,5,11,101,103,107,109,113,127,131,137,139.$$ 
\end{example}

\noindent
\textbf{Note added on August 11, 2010.} After our work was submitted, H. Brenner and A. Kaid  wrote the paper \cite{BK} providing, by means of a nice geometric argument, a more explicit characterization, in the special case $\alpha =\beta =\gamma =d$, of the primes $p$ of our Theorem \ref{main} for which the algebra $A=K[x,y,z]/(x^{\alpha},y^{\beta},z^{\gamma })$ has the WLP in characteristic $p$. In particular, Brenner and Kaid were able to solve (positively) our Conjecture \ref{char2}. It remains an open and very interesting problem to find a combinatorial proof of Conjecture \ref{char2}. We are grateful to Holger Brenner for sending us a copy of \cite{BK}.

Also, C. Chen (Berkeley), A. Guo (Duke), X. Jin (Minnesota) and G. Liu (Princeton), four REU students working in Summer 2010 at the University of Minnesota under the direction of Vic Reiner and Dennis Stanton, found, among other interesting things, a beautiful combinatorial explanation for our connection WLP-plane partitions, which was only proved algebraically in this paper. We are grateful to Vic Reiner for sending us an early copy of their work \cite{CGJL}.

\section*{Acknowledgements} We warmly thank Junzo Watanabe for several insightful comments, Uwe Nagel, Juan Migliore and Adam Van Tuyl for personal communications, and Erik Stokes, a postdoc of the second author at Michigan Tech, for technical assistance with a computer program the first author has created. This paper is the result of a senior undergrad project the first author has done at Michigan Tech in the academic year 2009-10. The second author wants to acknowledge to be indeed the second author --- the first author's overall contribution to this paper is above 50\%.\\



\begin{thebibliography}{ll}

\bibitem{And} G. Andrews: \lq \lq The theory of Partitions'', Encyclopedia of Mathematics and its Applications, Vol. II, Addison-Wesley, Reading, Mass.-London-Amsterdam (1976).

\bibitem{AE} G. Andrews and K. Eriksson: \lq \lq Integer Partitions'', Cambridge University Press, Cambridge, U.K. (2004).

\bibitem{bona} M. B\'ona: \lq \lq Introduction to enumerative combinatorics'', with a foreword by Richard Stanley. Walter Rudin Student Series in Advanced Mathematics, McGraw Hill Higher Education, Boston, MA (2007).

\bibitem{BK} H. Brenner and A. Kaid: \emph{A note on the weak Lefschetz property of monomial complete intersections in positive characteristic}, Collect. Math., to appear.

\bibitem{CGJL} C. Chen, A. Guo, X. Jin and G. Liu: \emph{Trivariate monomial complete intersections and plane partitions}, preprint. Available on the arXiv at \url{http://arxiv.org/abs/1008.1426}.

\bibitem{3} D. Cook II and U. Nagel: {\em The Weak Lefschetz Property, Monomial Ideals,
and Lozenges}, preprint. Available on the arXiv at \url{http://arxiv.org/abs/0909.3509v1}.

\bibitem{DT} G. David and C. Tomei: {\em The problem of the calissons}, Amer. Math. Monthly {\bf 96} (1989), 429-431.

\bibitem{DZZ} P. Di Francesco, P. Zinn-Justin and J.-B. Zuber: {\em A bijection between classes of fully packed loops and plane partitions}, Electron. J. Combin. {\bf 11} (2004), no. 1, Research Paper 64, 11 pp.

\bibitem{hawa} M. Hara and J. Watanabe: {\em The determinants of certain matrices arising from the Boolean lattice}, Discrete Math. {\bf 308} (2008), no. 23, 5815-5822.

\bibitem{HMNW} T. Harima, J. Migliore, U. Nagel and J. Watanabe: {\em The Weak and Strong Lefschetz Properties for Artinian $K$-Algebras}, J.  Algebra {\bf 262} (2003), 99-126.

\bibitem{Jacobi} D.E. Knuth: {\em Overlapping pfaffians}, Electron. J. Combin. \textbf{3} (1996), no. 2 (``The Foata Festschrift"), 11-24.

\bibitem{4} C. Krattenthaler: {\em Advanced determinant calculus}, S´eminaire Lotharingien Combin. \textbf{42} (1999), Article B42q, 67 pp.

\bibitem{kra} C. Krattenthaler: {\em Another involution principle-free bijective proof of Stanley's hook-content formula}, J. Combin. Theory Ser. A \textbf{88} (1999), no. 1, 66-92.

\bibitem{Macdonald} I.G. Macdonald: \lq \lq Symmetric Functions and Hall Polynomials'', Oxford University Press (1999).

\bibitem{Percy} P.A. MacMahon: \lq \lq Combinatorial analysis'', Chelsea Publishing Co., New York, NY (1960).

\bibitem{6} Mathematica, {\em A computational software program used in scientific, engineering, and mathematical fields and other areas of technical computing}. Available at \url{http://www.wolfram.com}.

\bibitem{5} J. Migliore, R. Mir\`o-Roig and U. Nagel: {\em Almost complete intersections and the Weak Lefschetz Property},  Trans. Amer. Math. Soc., to appear.

\bibitem{MZ} J. Migliore and F. Zanello: {\em The Hilbert functions which force the Weak Lefschetz Property}, J. Pure Appl. Algebra {\bf 210} (2007), no. 2, 465-471.

\bibitem{MS} E. Miller and B. Sturmfels: \lq \lq Combinatorial commutative algebra'', Graduate Texts in Mathematics, No. 227, Springer-Verlag, New York (2005).

\bibitem{pak} I. Pak: {\em Partition bijections, a survey},  Ramanujuan J. {\bf 12} (2006), 5-75.
 
\bibitem{RRR} L. Reid, L.G. Roberts and M. Roitman: {\em On Complete Intersections and their Hilbert Functions}, Canad. Math. Bull. \textbf{34} (1991), 525-535.

\bibitem{robe} P. Roberts: {\em A computation in local cohomology}, Contemp. Math. {\bf 159} (1994), 351-356.

\bibitem{St2} R. Stanley: {\em Hilbert functions of graded algebras}, Adv. Math. {\bf 28} (1978), 57-83.

\bibitem{St0} R. Stanley: \lq \lq Enumerative Combinatorics'', Vol. I, Second Ed., Cambridge Studies in Advanced Mathematics, Cambridge University Press (2010).

\bibitem{watanabe} J. Watanabe: {\em The Dilworth number of Artinian rings and finite posets with rank function}, Commutative Algebra and Combinatorics, Advanced Studies in Pure Math., Vol. 11, Kinokuniya Co. North Holland, Amsterdam (1987), 303-312.

\bibitem{ZZ} F. Zanello and J. Zylinski: {\em Forcing the Strong Lefschetz and the Maximal Rank Properties}, J. Pure Appl. Algebra {\bf 213} (2009), no. 6, 1026-1030.

\end{thebibliography}
\end{document}